\documentstyle[12pt]{article}
\newtheorem{theo}{Theorem}

\newtheorem{prop}[theo]{Proposition}
\def\beq{\begin{equation}}
\def\eeq{\end{equation}}
\def\bea{\begin{eqnarray}}
\def\eea{\end{eqnarray}}
\def\nn{\nonumber}
\def\ba{\begin{array}}
\def\ea{\end{array}}
\def\lb{[\![}
\def\rb{]\!]}

\def\ep{\epsilon}
\def\vep{\varepsilon}
\font\twelvemsbm=msbm10 at 12 true pt
\font\eightmsbm=msbm8
\font\sevenmsbm=msbm7
\newfam\msbmfam
\textfont\msbmfam=\twelvemsbm
\scriptfont\msbmfam=\eightmsbm
\scriptscriptfont\msbmfam=\sevenmsbm
\def\Bbb#1{{\fam\msbmfam\relax#1}} 
\def\C{\Bbb C}
\def\Zah{\Bbb Z}
\begin{document}
\renewcommand{\thefootnote}{\fnsymbol{footnote}}
\setcounter{footnote}{1}

\begin{center}
{\Large \bf The quantum superalgebra $U_q[osp(1/2n)]$~: deformed
para-Bose operators and root of unity representations}

\vspace{2cm}
{\bf T.D.\ Palev\footnote{Permanent address: Institute for
Nuclear Research and Nuclear Energy, 1784 Sofia, Bulgaria;
E-mail~: palev@bgearn.bitnet}
and
J.\ Van der Jeugt\footnote{Senior Research Associate of N.F.W.O.
(National Fund for Scientific Research of Belgium);
E-mail address: Joris.VanderJeugt@rug.ac.be}}

\vspace{0.5cm}
Department of Applied Mathematics and Computer Science \\
University of Ghent, Krijgslaan 281-S9, B-9000 Gent, Belgium
\end{center}

\vspace{2cm}

\noindent
{\bf Abstract}:
We recall the relation between the Lie superalgebra $osp(1/2n)$
and para-Bose operators. The quantum superalgebra
$U_q[osp(1/2n)]$, defined as usual in terms of its Chevalley
generators, is shown to be isomorphic to an associative algebra
generated by so-called pre-oscillator operators satisfying a
number of relations. From these relations, and the analogue with
the non-deformed case, one can interpret these pre-oscillator
operators as deformed para-Bose operators. Some consequences for
$U_q[osp(1/2n)]$ (Cartan-Weyl basis, Poincar\'e-Birkhoff-Witt
basis) and its Hopf subalgebra $U_q[gl(n)]$ are pointed out.
Finally, using a realization in terms of ``$q$-commuting''
$q$-bosons, we construct an irreducible finite-dimensional unitary
Fock representation of $U_q[osp(1/2n)]$ and its decomposition in
terms of $U_q[gl(n)]$ representations when $q$ is a root of unity.

\vspace{1cm}

\noindent
Short title: Quantum superalgebra $U_q[osp(1/2n)]$

\vspace{1cm}

\noindent
Phys. Abs. Class. Nos.~: 0210, 0220, 0365

\newpage

\renewcommand{\theequation}{\arabic{section}.{\arabic{equation}}}
\setcounter{equation}{0}
\section{Introduction}

It has been established in a number of papers
\cite{Palev78,Ohnuki,Palev82,Jag,Blank,Biswas,FlatoF,Okubo}
that the proper
mathematical framework of a system of $n$ para-Bose operators is
the theory of (unitarizable) representations of
the Lie superalgebra $osp(1/2n)$. This Lie superalgebra is
$B(0/n)$ in Kac's notation, and its finite-dimensional irreducible
representations (irreps) are completely classified (even their characters
are known). Unfortunately, these representations are not
unitarizable (see e.g.\ Ref.~\cite{Hughes}).
The infinite-dimensional unitarizable
representations have not been classified, and so far only
certain special cases (corresponding to parastatistics of order
$p$) have been considered.

Following the recent interest in $q$-deformations, a number of
papers have dealt with deformed parastatistics and in particular
with $q$-deformed para-Bose operators from different points of
view~\cite{Greenberg87}--\cite{Green94}. The definition
of $q$-deformed para-Bose operators depends highly upon the
framework one is working in, and most often one is inspired by
considering $q$-analogues of ordinary Bose operators.
Such approaches usually lead to deformed para-Bose operators that are
incompatible with a Hopf algebra structure.

{}From our point of view, the natural {\em ansatz} is the
equivalence between $osp(1/2n)$ irreps and para-Bose
operator representations. The $q$-deformed superalgebra
$U_q[osp(1/2n)]$, a Hopf superalgebra, is by now a classical
concept~\cite{Chaichian}--\cite{dHoker}.
Inside $U_q[osp(1/2n)]$ one can in a natural way
identify a set of elements with $q$-deformed para-Bose
operators. This leads to two basic results~: an alternative
definition of $U_q[osp(1/2n)]$ in terms of non-Chevalley
generators (the deformed para-Bose operators) satisfying a
number of relations; and the $q$-analogue of all triple
relations defining para-Bose statistics. These relations, on
their turn, imply the existence of a Poincar\'e-Birkhoff-Witt
theorem for $U_q[osp(1/2n)]$ and thus a basis in terms of
normally ordered monomials of Cartan-Weyl generators, the
expressions of which become extremely simple in terms of the
deformed para-Bose operators. Apart from these consequences, we
also state a number of results for the Hopf subalgebra
$U_q[gl(n)]\subset U_q[osp(1/2n)]$ and its realization in terms
of the deformed para-Bose operators.

Just as the defining relations for ordinary para-Bose operators
are automatically satisfied by canonical Bose operators (leading
to the oscillator representation of $osp(1/2n)$), we construct
in the present paper a set of deformed Bose operators that
satisfy, as a particular example, the deformed para-Bose
operator relations. These deformed Bose operators are (up to
some factor) equal to the usual $q$-bosons, but different modes
``$q$-commute'' instead of being commutative. The Fock space for
these operators is constructed, and shown to be unitarizable
only if $q$ is a primitive root of unity, in which case the
representation becomes finite-dimensional. These
finite-dimensional representations are given explicitly, and
their decomposition into irreducible representation of
$U_q[gl(n)]$ is considered.

\section{The para-Bose algebra $pB_n$ and its relation to the
Lie superalgebra $osp(1/2n)$}
\setcounter{equation}{0}

Let $A^\pm_i$ ($i=1,\ldots,n$) be a system of $n$ para-Bose
operators. The defining relation for para-Bose operators (parabosons),
introduced in quantum field theory by Green~\cite{Green} as a
possible generalization of integer spin field statistics
(see~\cite{Ohnuki} for a general introduction to
parastatistics), is given by~:
\beq
[ \{A_i^\xi,A_j^\eta \},A_k^\ep ]=
(\ep -\eta)\delta_{jk}A_i^\xi
+(\ep -\xi )\delta_{ik}A_j^\eta,\quad(\xi,\eta,\ep=\pm 1).
\label{pBose}
\eeq
The relations~(\ref{pBose}) generalize the canonical commutation
relations of ordinary Bose operators (bosons) $a^\pm_i$~:
\beq
[a_i^-,a_j^+]=\delta_{ij}, \quad [a_i^-,a_j^-]=[a_i^+,a_j^+]=0,
\label{Bose}
\eeq
and it is trivial to verify that the $a_i^\pm$ do indeed
satisfy~(\ref{pBose}).

The para-Bose algebra $pB_n$ is defined as the associative
algebra with unity over $\C$ with generators $A^\pm_i$, subject to the
relations~(\ref{pBose}). In fact, $pB_n$ is turned into a superalgebra
(associative $\Zah_2$-graded algebra) by the requirement
$\deg(A^\pm_i) = \bar 1$, $\forall i\in\{1,\ldots,n\}$, where
$\Zah_2 =\{\bar 0, \bar 1\}$. By defining the supercommutator
between any two homogeneous elements $a$ and $b$ of $pB_n$ by
\beq
\lb a,b \rb =ab-(-1)^{deg(a)deg(b)}ba ,
\label{supercomm}
\eeq
and extending it by bilinearity to the whole algebra, $pB_n$ is
turned into a Lie superalgebra. Thus $n$ pairs of para-Bose
operators generate a Lie superalgebra~\cite{Omote}, and we shall
recall that in the present case this Lie superalgebra can be
identified with $osp(1/2n)$~\cite{Palev78}.

For this purpose, define the Lie superalgebra $osp(1/2n)$ (or
$B(0/n)$ in Kac's notation~\cite{Kac}) as the set of
$(2n+1)\times (2n+1)$ complex matrices of the form
\beq
\left(\begin{array}{ccc}
          0     & x & y \\
	  y^T & d & e   \\
	 -x^T   & f &-d^T
\end{array}\right),
\label{matrix}
\eeq
where $T$ stands for transposition, $d,e$ and $f$ are $n\times
n$ matrices with $e^T=e$ and $f^T=f$, and $x,y$ are $1\times n$
matrices.
The even subalgebra $osp(1/2n)_{\bar 0}$ consists of all matrices with
$x=y=0$ and is isomorphic to the symplectic Lie algebra
$sp(2n)$.  The odd subspace $osp(1/2n)_{\bar 1}$ consists of all
matrices with $d=e=f=0$. The supercommutator between homogeneous
elements is defined by means of~(\ref{supercomm}), and extended
by bilinearity. Let $E_{kl}$ denote the $(2n+1)\times(2n+1)$ matrix
with 1 at the intersection of row $k$ and column $l$, and
zero elsewhere (where rows and columns are labelled from 0 to
$2n$). The Cartan subalgebra $H$ of $osp(1/2n)$ is spanned by the
elements $H_i=-E_{ii}+E_{n+i,n+i}$ ($i=1,\ldots,n$). With a
suitable basis $\vep_i$ ($i=1,\ldots,n$) of the dual space $H^*$,
the roots of $osp(1/2n)$ consist of $2n$ odd roots $\pm\vep_i$
and $2n^2$ even roots $\pm\vep_i\pm\vep_j$~\cite{Kac}. The Lie
superalgebra $osp(1/2n)$ has the usual root space decomposition,
with all root spaces one-dimensional. The root vectors
corresponding to the odd roots can be written as follows~:
\bea
\vep_i \quad&:&\quad A_i^- =\sqrt{2} (E_{0,i}-E_{i+n,0}),\quad
i=1,\ldots,n, \label{A-}\\
-\vep_i\quad&:&\quad A_i^+ =\sqrt{2} (E_{0,i+n}+E_{i,0}), \quad
i=1,\ldots,n. \label{A+}
\eea
Therefore the anticommutator $\{A_i^{-\xi},
A_j^{-\eta}\}$ is a root vector with root $\xi\vep_i+\eta\vep_j$
($\xi,\eta=\pm$), and one finds $\{A^-_i,A^+_i\}=-2H_i$
($i,j=1,\ldots,n$). This implies that the odd root vectors
$A^\pm_i$ generate the whole Lie superalgebra $osp(1/2n)$.
Moreover, one can verify that the elements (\ref{A-})-(\ref{A+})
satisfy indeed~(\ref{pBose}). This gives rise to the following

\begin{prop}
The para-Bose algebra $pB_n$ is isomorphic to the universal
enveloping algebra $U[osp(1/2n)]$ of $osp(1/2n)$. The finite
dimensional subspace
\beq
{\rm span } \{\{A_i^\xi,A_j^\eta \}, A_k^\ep
 \vert  i,j,k=1,\ldots,n;\; \xi,\eta,\ep=\pm\}
\label{osp12n},
\eeq
endowed with the supercommutator~(\ref{supercomm}), is a
Lie superalgebra isomorphic to $osp(1/2n)$.
\end{prop}

Note that the even elements of $osp(1/2n)$ are spanned by the
$\{A_i^\xi, A_j^\eta\}$, thus $sp(2n)={\rm span}\{ \{A_i^\xi,
A_j^\eta\} | i, j = 1, \ldots,n;\; \xi,\eta=\pm\}$. In particular,
{}from~(\ref{pBose}) one can derive a compact expression for the
commutation relations of the $sp(2n)$ basis elements~:
\bea
[\{A_i^\xi,A_j^\eta\},\{A_k^\epsilon,A_l^\varphi\} ]& =&
(\epsilon -\eta)\delta_{jk}\{A_i^\xi,A_l^\varphi\}+
(\epsilon -\xi )\delta_{ik}\{A_j^\eta,A_l^\varphi\} \nn\\
& &+(\varphi-\eta)\delta_{jl}\{A_i^\xi,A_k^\epsilon\}+
(\varphi-\xi)\delta_{il}\{A_j^\eta,A_k^\epsilon\}.
\label{sp2n}
\eea

As a consequence from Proposition~1, determining all representations of the
para-Bose operators is completely equivalent with finding all
representations of the Lie superalgebra $osp(1/2n)$. For
finite-dimensional irreducible representations of $osp(1/2n)$,
there exists a character formula~\cite{Kac}, but explicit
formulae for matrix elements are not available in the
literature. In para-Bose statistics, one is rather
interested in infinite dimensional unitarizable representations.
The Fock space corresponding to the ordinary Bose
operators~(\ref{Bose}) provides one such example; representations
corresponding to a fixed order of
parastatistics~\cite{Green,Greenberg} provide in principle others; but apart
{}from that no general theory of unitarizable representations of
$osp(1/2n)$ exists.

Besides the definition of $U[osp(1/2n)]$ in terms of generators
$A^\pm_i$ subject to the relations~(\ref{pBose}), there is an
alternative definition in terms of Chevalley generators. This is
perhaps the definition that most readers are more familiar with.
Let $(\alpha_{ij})$ be a Cartan matrix chosen
as an $n  \times n $ symmetric matrix with
\beq
\alpha_{nn}=1,\quad \alpha_{ii}=2, \quad
\alpha_{i,i+1}=\alpha_{i+1,i}=-1, \quad i=1,\ldots,n-1,\quad
{\rm all \; other}\; \alpha_{ij}=0.
\label{cartan}
\eeq
Again $U[osp(1/2n)]$ is defined as an associative superalgebra
in terms of a number of generators subject to relations. The
generators are the elements $h_i, e_i, f_i$ ($i=1,\ldots,n$);
the relations are the Cartan-Kac relations
\bea
&&[h_i,h_j]=0,\nn \\
&&[h_i,e_j]=\alpha_{ij}e_j,\quad
  [h_i,f_j]=-\alpha_{ij}f_j,\quad \nn\\
&&[e_i,f_j]=\delta_{ij}h_i\quad{\rm except\ for}\quad i=j=n, \label{CK}\\
&& \{ e_n,f_n\}=h_n ; \nn
\eea
the $e$-Serre relations
\bea
&& [e_i,e_j]=0, \quad{\rm for}\quad |i-j|>1,\nn\\
&& e_i^2e_{i+1}-2e_ie_{i+1}e_i+e_{i+1}e_i^2=0,
   \quad i=1,\ldots,n-1, \nn\\
&& e_i^2e_{i-1}-2e_ie_{i-1}e_i+e_{i-1}e_i^2=0,
   \quad i=2,\ldots,n-1,\label{eSerre} \\
&& e_n^3e_{n-1}- (e_n^2e_{n-1}e_n+e_ne_{n-1}e_n^2)+
  e_{n-1}e_n^3=0; \nn
\eea
and the $f$-Serre relations obtained from (\ref{eSerre}) by replacing
everywhere $e_i$ by $f_i$. The grading on the superalgebra is
induced from the grading on its generators~: $\deg (e_n)=\deg(f_n)=\bar 1$
and $\deg(e_i)=\deg(f_i)=\bar 0$ for $i=1,\ldots,n-1$.
Although the present definition of $U[osp(1/2n)]$ looks somewhat
more complicated than the one given in terms of the para-Bose
operators $A_i^\pm$, it cannot be avoided when turning to a Hopf
superalgebra deformation $U_q[osp(1/2n)]$. Such deformations are
well defined only in a Chevalley basis. In fact, relations like
the ones given below, expressing the para-Bose operators in
terms of the Chevalley generators, will be reconsidered in the
following section to yield a definition of deformed para-Bose
operators in terms of the Chevalley generators of
$U_q[osp(1/2n)]$. Here, the relations read~:
\bea
&&A_i^-=-\sqrt{2}
[e_i,[e_{i+1},[e_{i+2},[\ldots,[e_{n-2},[e_{n-1},e_n]
\ldots ],\quad i=1,\ldots,n-1 \nn\\
&&A_i^+=\sqrt{2}
  [\ldots [f_n,f_{n-1}]_{q},f_{n-2}],\ldots],
  f_{i+2}],f_{i+1}],f_{i}], \quad i=1,\ldots,n-1, \\
&&A_n^-=-\sqrt{2}e_n, \qquad A_n^+=\sqrt{2}f_n. \nn
\eea
Clearly, the expressions of the Chevalley generators in terms of
para-Bose operators also exist~:
\bea
&&e_i={1\over 2}\{A_i^-,A_{i+1}^+\},\quad
f_i={1\over 2}\{A_i^+,A_{i+1}^-\},\quad i=1,\ldots,n-1,\label{ei}\\
&&h_i={1\over 2}\{A_{i+1}^-,A_{i+1}^+\}-
  \{A_{i}^-,A_{i}^+\}, \quad i=1,\ldots,n-1, \label{hi}\\
&&e_n=-{1\over \sqrt 2}A_n^-,
  \quad f_n={1\over \sqrt 2}A_n^+,
  \quad h_n=-{1\over 2}\{A_n^-,A_n^+\}.
\eea
An important statement can be deduced from these relations. From
(\ref{CK})-(\ref{eSerre}) it follows that the
enveloping algebra of $gl(n)$ is generated by $h_i$
($i=1,\ldots,n$) and $e_i,f_i$ ($i=1,\ldots,n-1$). Then
(\ref{ei})-(\ref{hi}) and (\ref{sp2n}) show that $U[gl(n)]$, a subalgebra of
$U[osp(1/2n)]$, is generated by the elements $\{A_i^-,A_j^+\}$
($i,j=1,\ldots,n$). In other words~:
\beq
gl(n)={\rm span }\{\{A_i^-,A_j^+\}| i,j=1,\ldots,n\},
\label{gln}
\eeq
endowed with the usual commutator product. Replacing in
(\ref{gln}) the para-Bose operators with the ordinary Bose
operators (\ref{Bose}) yields the familiar Schwinger realization
of $gl(n)$; but observe that in the superalgebra grading
the para-Bose and Bose operators are odd (``fermionic'')
operators rather than even (``bosonic'') operators.

\setcounter{equation}{0}
\section{The quantum superalgebra $U_q[osp(1/2n)]$ and
$q$-de\-for\-med para-Bose operators}

In the present section we consider the well known quantum
superalgebra $U_q[osp(1/2n)]$. A general procedure to construct
a so-called Cartan-Weyl basis of $U_q[osp(1/2n)]$ (the analogue
of a Cartan-Weyl basis for the Lie superalgebra $osp(1/2n)$) was
given in Ref.~\cite{Khoroshkin}. We follow this procedure, and
identify -- as in the non-deformed case of Section~1 -- the
Cartan-Weyl basis elements corresponding to the odd roots as
``deformed para-Bose operators.'' The remaining questions that
are treated in this section are~: the alternative definition of
$U_q[osp(1/2n)]$ in terms of deformed para-Bose operators; the
triple relations between deformed para-Bose operators (the
analogue of~(\ref{pBose})); a Poincar\'e-Birkhoff-Witt (PBW) theorem for
$U_q[osp(1/2n)]$; the Cartan-Weyl basis of the subalgebra
$U_q[gl(n)]$ and its consequences.

First, we introduce $U_{q}\equiv U_{q}[osp(1/2n)]$ by means
of its classical definition in terms of the Cartan
matrix~(\ref{cartan}). $U_{q}$ is the associative superalgebra
with unity
over $\C$, generated by the elements $k_i^\pm=q^{\pm h_i}$,
$e_i, f_i$ ($i=1,\ldots,n$), subject to the following
relations~: the Cartan-Kac relations
\bea
&&k_ik_i^{-1}=k_i^{-1}k_i=1, \quad k_ik_j=k_jk_i, \nn \\
&&k_ie_j=q^{\alpha_{ij}}e_jk_i, \quad
k_if_j=q^{-\alpha_{ij}}f_jk_i, \nn\\
&&[e_i,f_j]=\delta_{ij}(k_i-k_i^{-1})/(q-q^{-1})
\quad {\rm except\ for}\quad i=j=n, \label{qCK}\\
&& \{ e_n,f_n\} =(k_n-k_n^{-1})/(q-q^{-1}); \nn
\eea
the $e$-Serre relations
\bea
&&[e_i,e_j]=0, \quad {\rm for}\quad |i-j|>1,  \nn\\
&&e_i^2e_{i+1}-(q+q^{-1})e_ie_{i+1}e_i+e_{i+1}e_i^2=0,
   \quad i=1,\ldots,n-1, \nn\\
&&e_i^2e_{i-1}-(q  +q^{-1})e_ie_{i-1}e_i+e_{i-1}e_i^2=0,
   \quad i=2,\ldots,n-1, \label{qeSerre} \\
&&e_n^3e_{n-1}+(1-q  -q^{-1})(e_n^2e_{n-1}e_n+e_ne_{n-1}e_n^2)+
e_{n-1}e_n^3=0;\nn
\eea
and the $f$-Serre relations obtained from above by replacing
everywhere $e_i$ by $f_i$.
The grading on $U_{q}$ is induced from the requirement that the
generators $e_n,\; f_n$ are odd and all other generators are
even. It is known that $U_{q}$ can be endowed with a
comultiplication $\Delta$, a counit $\vep$ and an antipode $S$,
turning it into a
Hopf superalgebra; here we shall not be concerned with this
additional structure.

Following the procedure outlined in Ref.~\cite{Khoroshkin} one
determines the Cartan-Weyl elements corresponding to the odd
roots, and thus we define the deformed para-Bose operators as
follows (see also~\cite{Palev93A})~:
\bea
&& A_i^-=-\sqrt{2}
[e_i,[e_{i+1},[e_{i+2},[\ldots,[e_{n-2},[e_{n-1},e_n]_{q^{-1}}
]_{q^{-1}}\ldots ]_{q^{-1}},\nn\\
&& A_i^+=\sqrt{2}
  [\ldots [f_n,f_{n-1}]_{q},f_{n-2}]_{q},\ldots]_{q},
f_{i+2}]_{q},f_{i+1}]_{q},f_{i}]_{q}, \quad
i=1,\ldots,n-1, \label{qAe}\\
&& A_n^-=-\sqrt{2}E_n, \quad
 A_n^+=\sqrt{2}F_n, \nn
\eea
where $[u,v]_q=uv-qvu$. Besides these, we also
introduce $n$ even ``Cartan" elements
\beq
L_i=k_ik_{i+1}\ldots k_n=q^{H_i} {\rm \ where\ }
H_i=h_i+h_{i+1}+\ldots +h_n, \quad i=1,\ldots,n.
\label{Li}
\eeq
We shall call the set of operators $A_i^\pm, L_i^\pm$
pre-oscillator operators, for reasons that will be obvious in
the following section. Using their definition
(\ref{qAe})-(\ref{Li}) and the defining relations
(\ref{qCK})-(\ref{qeSerre}), one can also determine the
Chevalley generators in terms of the pre-oscillator operators
($i\ne n$)~:
\bea
&&  e_i=-{q\over 2}\{A_i^-,A_{i+1}^+\}L_{i+1}^{-1},\quad
  f_i=-{1\over 2q}L_{i+1}\{A_i^+,A_{i+1}^-\}, \nn\\
&& e_n=-(2)^{-1/2}A_n^-,\quad
 f_n=(2)^{-1/2}A_n^+ . \label{qeA}
\eea

Now, one can determine the relations between the pre-oscillator
operators. A set of such relations were already obtained in
Ref.~\cite{Palev93A}, but with the purpose of giving an alternative
definition of $U_{q}[osp(1/2n)]$ in terms of deformed para-Bose
or pre-oscillator operators it would be interesting to find the
direct analogue of (\ref{qCK})-(\ref{qeSerre}), and thus present a
minimal set of relations. This is given in the following~:

\begin{prop}
The relations of  $U_{q}[osp(1/2n)]$
in terms of its Chevalley generators
hold if and only if the pre-oscillator operators satisfy
($i,j=1,\ldots,n$, $\xi=\pm$)~:
\bea
&& L_iL_i^{-1}=L_i^{-1}L_i=1, \quad L_iL_j=L_jL_i,  \nn\\
&& L_iA_j^{\pm}={q}^{\mp\delta_{ij}}A_j^{\pm}L_i, \nn\\
&& \{A_i^-,A_i^+\}=-2(L_i- L_i^{-1})/(q-q^{-1}), \label{minim}\\
&& [\{A_i^{-\xi},A_{i\pm 1}^{\xi}\},A_j^{-\xi}]_{{q}^{\pm
  \delta _{ij}}}=-2\xi \delta _{j,i \pm 1}L_j^{\pm \xi}A_i^{-\xi},\nn\\
&& [\{A_{n-1}^\xi,A_n^\xi\},A_n^\xi]]_{q}=0 . \nn
\eea
\end{prop}

We shall refer to (\ref{minim}) as the pre-oscillator
realization of $U_{q}$.
Although the definition of $U_{q}[osp(1/2n)]$ in terms of deformed
para-Bose operators is more complicated than in the
non-deformed case, where only one relation~(\ref{pBose}) was
necessary, it should be observed that the present definition
by means of (\ref{minim}) is certainly not more involved than
the list of classical relations, i.e.\ the Cartan-Kac, the
$e$-Serre and the $f$-Serre relations. The remaining advantage
of the definition by means of Chevalley generators is the
simplicity of the other Hopf superalgebra functions $\Delta,
\vep$ and $S$ which become very complicated expressions on
$A_i^\pm$ and $L_i^{\pm 1}$, although they are certainly well defined.

Next, we are concerned with deriving the analogue of
(\ref{pBose}), i.e.\ all triple relations between deformed
para-Bose operators. These relations are derived using
(\ref{minim}), and they are far more complicated than the
classical relation (\ref{pBose}). In fact, it would be quite
impossible to cast all of them in a single expression.
Nevertheless, we think they are of importance since they
constitute the direct Hopf algebra $q$-deformation of para-Bose statistics.
In the following list, we use the abbreviation~:
\beq
\tau_{i_1,i_2,\ldots,i_k}=
\left\{ \begin{array}{rl}
        -1, & {\rm if\ }i_1>i_2>\ldots >i_k,\\
        1, & {\rm if\ }i_1<i_2<\ldots <i_k,\\
	0,  & {\rm otherwise.}
        \end{array}\right.
\eeq
Then, this list reads~:
\bea
&& [\{A_i^{-\xi},A_j^\xi\},A_k^\xi]_{q^{\tau_{ji}\delta_{jk}}} =
2\xi\delta_{ik}A_j^{\xi}L_k^{\xi\tau_{ij}} -
(q-q^{-1})\tau_{ikj} \{A_i^{-\xi},A_k^\xi\}A_j^\xi,\quad(i\ne
j), \label{1}\\
&&
[\{A_i^\xi,A_j^\xi\},A_k^{-\xi}]=
(q-q^{-1})\left(\tau_{kij}\{A_i^\xi,A_k^{-\xi}\}A_j^\xi +
\tau_{kji}\{A_j^\xi,A_k^{-\xi}\}A_i^\xi \right) \nn\\
&&\qquad-2\xi\delta_{ik}A_j^\xi L_i^{\xi\tau_{ij}}
  -2\xi\delta_{jk}A_i^\xi L_j^{\xi\tau_{ji}}, \quad(i\ne j),\label{2}\\
&&[\{A_i^\xi,A_i^\eta\},A_k^{-\eta}]=
2\delta_{\xi\eta}(q^{\tau_{ki}}-1)\{A_i^\xi,A_k^{-\xi}\}
A_i^\xi \nn\\
&& \qquad -2\xi\eta(1+\delta_{\xi\eta})\delta_{ik}
A_i^\xi\left( (q^{\xi}-1) L_i^{-1}+(1-q^{-\xi})L_i\right)/
(q-q^{-1}) , \label{3}\\
&&[\{A_i^\xi,A_j^\xi\},A_k^\xi]_{q^{\tau_{ik}+\tau_{jk}}}=0. \label{4}
\eea

Apart from giving the $q$-analogue of para-Bose relations,
(\ref{1})-(\ref{4}) also imply that monomials in $A_i^\pm$ and
$\{A_k^\xi, A_l^\eta\}$ can be reordered. A detailed
investigation of the quadruple relations, i.e.\ the relations
between $\{A_i^\xi, A_j^\eta\}$ and $\{A_k^\ep, A_l^\zeta\}$
(which shall not be given here, since the complete list is too
long), yields the following
\begin{prop}
The set of operators
\beq
L_i^{\pm 1}, \;A_i^\pm,\;
\{A_i^-,A_j^+ \},\; \{A_k^\xi,A_l^\xi \},\quad i\neq j,\quad
i,j,k,l=1,\ldots,n,\quad \xi=\pm,
\label{CW}
\eeq
give a Cartan-Weyl basis of $U_{q}[osp(1/2n)]$.
The set of all normally ordered monomials~\cite{Khoroshkin}
in (\ref{CW}) constitute a basis in
$U_{q}[osp(1/2n)]$ (PBW theorem).
\label{proposp}
\end{prop}
This shows one of the important advantages of the deformed
para-Bose operators~: they yield a very simple basis for
$U_{q}[osp(1/2n)]$. At the same time we can restrict the above
statements to the subalgebra $U_{q}[gl(n)]$~:
\begin{prop}
The operators
\beq
L_i^{\pm 1}, \;
\{A_i^-,A_j^+ \} ,\quad i\neq j
=1,\ldots,n
\label{CWgln}
\eeq
constitute a Cartan-Weyl basis for the Hopf
superalgebra $U_{q}[gl(n)]$; the normally ordered monomials in
(\ref{CWgln}) form a basis in $U_{q}[gl(n)]$.
\label{propgln}
\end{prop}

Let us say a few words about the normal order
for the elements (\ref{CW}).
Usually, one takes an order $<$ such
that for positive root vectors ({\em prv}) $A_i^-$,
$\{A_i^-,A_j^\pm\}$ ($i<j$), for negative root vectors ({\em
nrv}), and for the Cartan generators $L_i$ the inequality ${\em
prv} < {\em nrv} < L_i$ holds. Among the {\em prv} the order is
taken to be~\cite{Palev93A}~:
\beq
\{A_i^-,A_k^\xi\}<\{A_i^-,A_l^\xi\},\quad
{\rm for}\quad k<l, \label{ord1}
\eeq
\beq
\{A_i^-,A_k^+\}<A_i^-<\{A_i^-,A_l^-\}
 <\{A_j^-,A_r^+\}<A_j^-<\{A_j^-,A_s^-\},
 \quad {\rm for}\quad i<j;\; i<k,l;\; j<r,s. \label{ord2}
\eeq
For the proof of Proposition~\ref{proposp} one has to show that
the unordered product of any two Cartan-Weyl elements (\ref{CW})
can be represented as a linear combination of normally ordered
products, and that this procedure of ordering is finite
when applied to a finite unordered monomial in the
Cartan-Weyl elements.
As we mentioned previously, this can be deduced from
the triple relations (\ref{1})-(\ref{4}), and from a list
of quadruple relations which is too long to be included here.
When restricting the elements to (\ref{CWgln}), it turns out
that also the quadruple relations can be summarized rather
easily, yielding thus a complete proof of
Proposition~\ref{propgln}.

For this purpose, introduce the function
\beq
\theta_{i_1,i_2,\ldots,i_k}=\left\{
 \begin{array}{rl}
  1 & {\rm if\ } i_1>i_2>\ldots >i_k,\\
  0 & {\rm otherwise.}
 \end{array}\right.
\eeq
A set of Cartan-Weyl elements of $U_{q}[gl(n)]$ has been
considered before~\cite{Palev91}, and consists of $n$ ``Cartan"
elements $L_1, L_2,\ldots,L_n$ and $n(n-1)$ root vectors $e_{ij}$,
$i\neq j=1,\ldots,n$. Remember that $e_{ij}$ is positive
if $i<j$ and negative if $i>j$.
Among the {\em prv}, the normal order induced
{}from~(\ref{ord1})-(\ref{ord2}) yields~:
\beq
e_{ij}<e_{kl},\;\; {\rm if}\;\; i<k \;\;{\rm or}\;\; i=k \;\;
{\rm and}\;\; j<l;
\label{ord3}
\eeq
for {\em nrv} one takes the same rule~(\ref{ord3}), and one
chooses ${\em prv}<{\em nrv}<L_i$. This yields a normal order
for the $U_{q}[gl(n)]$ Cartan-Weyl basis elements. A complete set
of relations is given by~:
\begin{enumerate}
\item
\beq
L_i^\xi L_j^\eta=L_j^\eta L_i^\xi, \quad
L_ie_{jk}=q^{\delta_{ij}-\delta_{ik}}e_{jk}L_i ; \label{g1}
\eeq
\item
For any $e_{ij}>0$  and $e_{kl}<0$~:
\bea
[e_{ij},e_{kl}]&=&\left( (q-q^{-1}) \theta_{jkil}e_{kj}e_{il}
 -\delta_{il}\theta_{jk}e_{kj}+
  \delta_{jk}\theta_{il}e_{il}\right) L_k L_i^{-1}  \nn\\
&&+L_l L_j^{-1}\left(-(q-q^{-1})\theta_{kjli}e_{il}e_{kj}
-\delta_{il}\theta_{kj}e_{kj}+\delta_{jk}\theta_{li}e_{il}
\right) \nn\\
&& + \delta_{il}\delta_{jk}
(L_i L_j^{-1} - L_i^{-1} L_j)/(q-q^{-1}); \label{g2}
\eea
\item
Set $\xi=1$ if  $0<e_{ij}<e_{kl}$, and $\xi=-1$ if
   $0>e_{ij}>e_{kl}$. Then
\beq
e_{ij}e_{kl}-q^{\xi(\delta_{ik}-\delta_{il}-\delta_{jk}+\delta_{jl})}
e_{kl}e_{ij}
=\delta_{jk}e_{il}+(q-q^{-1})\tau_{ljki}e_{kj}e_{il}.
\label{g3}
\eeq
\end{enumerate}

To obtain the above relations we have considered Equations
(3.10)-(3.14) from Ref.~\cite{Palev91} only for the even
generators, first replacing
$q$ by $-q$ and then $q^{e_{ii}}$ by $L_i$.
The link with the present operators is now given by~:
\beq
e_{ij}=-{1\over 2}L_j^{-1} \{A_i^-,A_j^+\} {\rm \ for\ } i<j,
 \quad {\rm and\ \ }
e_{ij}=-{1\over 2} \{A_i^-,A_j^+\}L_i{\rm \ for\ }
i>j.
\label{eij}
\eeq
It is tedious but straightforward to verify that (\ref{eij})
satisfy indeed the relations (\ref{g1})-(\ref{g3}).

To conclude this section, observe that for
$U_{q}[osp(1/2n)]$, resp.\ for $U_{q}[gl(n)]$, the expressions of the
Cartan-Weyl root vectors in terms of the deformed para-Bose
operators are the same as the expressions of the Cartan-Weyl
basis elements in terms of the non-deformed para-Bose operators
for $osp(1/2n)$, resp.\ for $gl(n)$.

\setcounter{equation}{0}
\section{Realization of $U_{q}[osp(1/2n)]$ in terms of deformed
Bose oscillators}

We have seen that the canonical Bose operators~(\ref{Bose})
satisfy the para-Bose relations~(\ref{pBose}). As a consequence
the Fock space built on the Bose operators forms a
representation of the para-Bose algebra $pB_n\equiv U[osp(1/2n)]$,
usually referred to as the oscillator representation of
$osp(1/2n)$. In the present section we shall define a new set of
operators which can be interpreted as deformed Bose operators.
They satisfy the deformed para-Bose relations~(\ref{minim}) and
hence their Fock space forms a representation of
$U_{q}[osp(1/2n)]$, the analogue of the usual oscillator
representation. Surprisingly, the unitarity conditions lead to
the condition that $q$ should be a root of unity, yielding
a finite-dimensional Fock space. This Fock representation is
studied in detail, and in particular its decomposition with
respect to $U_{q}[gl(n)]$ is considered, giving rise to certain
root of unity representations for this quantum algebra.

Define a  set of operators
\beq
a_i^\pm, \; \kappa_i=q^{N_i},\quad i=1,\ldots,n,
\eeq
satisfying the relations
\bea
&& a_i^-a_i^+ -q^{\pm 1}a_i^+a_i^-
  ={2\over {q^{1/2}+q^{-1/2}}}\kappa_i^{\mp 1}, \nn\\
&& \kappa_i a^{\pm}_j = q^{\pm\delta_{ij}} a^\pm_j \kappa_i, \label{a-t}\\
&& a_i^\xi a_j^\eta =q^{\xi \eta}
  a_j^\eta a_i^\xi , \quad {\rm for\ all\ }i<j. \nn
\eea
It is obvious that in the limit $q \rightarrow 1$ the operators
$a_i^\pm$ reduce to the usual Bose creation and annihilation
operators~(\ref{Bose}). For a fixed $i$ the relations coincide
(up to a multiple) with the usual so-called $q$-deformed
oscillators~\cite{Macfarlane,Biedenharn,Sun,Hayashi}. Note
however that the third relation in~(\ref{a-t}) implies that
different modes do not commute, but ``$q$-commute''; also such a
phenomenon has been considered
before~\cite{Pusz1,Pusz2,Hadjiivanov,Jagannathan,VdJ}.

Denote by $W_q(n)$  the associative algebra with unity over $\C$ with
generators $a_i^\pm,\; \kappa_i^{\pm 1}$ ($i=1,\ldots,n$) and
relations~(\ref{a-t}). Clearly, $W_q(n)$ is a deformation
of the canonical Weyl algebra $W(n)$, generated by $n$ pairs
of Bose operators.

\begin{prop}
The linear map $\varphi$ from $U_{q}[osp(1/2n)]$ into
$W_q(n)$, defined on the pre-oscil\-la\-tor generators as
\beq
\varphi( A_i^\pm) = a_i^\pm \quad
\varphi(L_i)=q^{-1/2}\kappa_i^{-1}\equiv l_i,
\quad i=1,\ldots,n,
\eeq
and extended on all elements by associativity is
a (associative algebra) homomorphism of $U_{q}[osp(1/2n)]$
onto $W_q(n)$.
\label{prop-al}
\end{prop}

The proof follows from the observation that the operators
$a_i^\pm$ and $l_i=q^{-1/2}\kappa_i^{-1}$ satisfy equations~(\ref{minim}).
Moreover the generators of $W_q(n)$ are among the images of
$\varphi$. In particular $\kappa_i=\varphi(q^{-1/2}L_i^{-1}) $. From
Propositions~\ref{proposp} and~\ref{prop-al} it follows that the
following elements yield an oscillator realization of the
Cartan-Weyl generators of $U_{q}[osp(1/2n)]$~:
\beq
l_i^{\pm 1}, \;a_i^\pm,\;
\{a_i^-,a_j^+ \},\; \{a_i^\xi,a_j^\xi \},\quad i\neq j
=1,\ldots,n.
\label{oscreal}
\eeq
Deformed Bose creation and annihilation operators have
been studied recently
\cite{Macfarlane}--\cite{Flato} and in the
past~\cite{Coon,Arik,Kuryshkin,Jannussis}. The presently
introduced operators~(\ref{a-t}) coincide with others only in
one mode. The main reason for introducing them lies in the
underlying connection with the deformed para-Bose
operators~(\ref{minim}). As an important consequence, we shall
see that unitary deformed oscillator representations exist only
for $q$ being a root of unity.

A representation of $U_{q}[osp(1/2n)]$ is said to be unitary if the
representation space is a Hilbert space and the representatives
of $A_i^\pm$ and $L_i^{\pm 1}$ (where $L_i^{\pm 1}=q^{\pm H_i}$ are supposed
to be diagonal) satisfy
\beq
(A_i^+)^\dagger =A_i^-,\;\; (H_i)^\dagger =H_i,
\label{herm}
\eeq
where $A^\dagger$ is the Hermitian conjugate to the operator $A$.
In particular, let us now consider the Fock space defined by
means of the deformed oscillators~(\ref{a-t}). Requiring
herein~(\ref{herm}) leads to $(a_i^+)^\dagger =a_i^-$,
$(N_i)^\dagger =N_i$. But the relations~(\ref{a-t}) remain
invariant under this conjugation if and only if $|q|=1$, i.e.\
if $q$ is a phase.

Next, we proceed with the construction of the Fock space. As
usual, the vacuum vector $|0\rangle$ is defined by means of
\beq
a_i^- |0\rangle=0 . \label{vac}
\eeq
The basis states are of the form
\beq
|m_1,\ldots,m_n\rangle = {\cal N}_{m_1,\ldots,m_n}
(a_1^+)^{m_1}\cdots (a_n^+)^{m_n} |0\rangle ,
\label{state}
\eeq
where ${\cal N}_{m_1,\ldots,m_n}$ is a normalization constant. Under
the unitarity condition, a simple calculation using~(\ref{a-t})
and~(\ref{vac}) leads to
\beq
\langle m_1,\ldots,m_n | m_1,\ldots, m_n\rangle =
|{\cal N}_{m_1,\ldots,m_n}|^2 \alpha(m_1)\cdots\alpha(m_n),
\label{norm}
\eeq
where
\beq
\alpha(m)={(q^{1/2}+q^{-1/2})^{m}\over{2^{m}}[m]_{q}!},
\label{alpha}
\eeq
and, as usual, $[x]_q=(q^x-q^{-x})/(q-q^{-1})$ and $[x]_q! =
[x]_q [x-1]_q \cdots [1]_q$. For all allowed values of the
labels $m_i$, the norm~(\ref{norm}) should be positive. In
particular, this implies that $\alpha(m_i)$ must be positive for
all $m_i=0,1,\ldots$. Then~(\ref{alpha}) implies that
$[m]_{q}$ should be positive for $m=0,1,\ldots$. But since
$q=e^{i\phi}$ is a pure phase,
$[m]_{q}=\sin(m\phi)/\sin(\phi)$. It follows that the only
admissible situation is when $q$ is a primitive root of unity,
\beq
q=e^{i\pi/k},
\eeq
for some positive integer $k$. In that case, the Fock space is
finite dimensional, with
\beq
m_i \in \{0,1,\ldots,k-1\}, \quad i=1,\ldots,n.
\eeq
Thus, the total dimension of the Fock space is $k^n$. When the
basis is taken to be orthonormal, i.e.\
\beq
{\cal N}_{m_1,\ldots,m_n}=\left( \alpha(m_1)\cdots\alpha(m_n)
\right)^{-1/2},
\eeq
the explicit action of the deformed oscillators, using the
short-hand notation
\bea
&& |m\rangle = | \ldots,m_{i-1},m_i,m_{i+1},\ldots\rangle \nn\\
&&
|m_i+1\rangle =|\ldots,m_{i-1},m_i+1,m_{i+1},\ldots\rangle,\label{shn} \\
&&
|m_i-1\rangle= |\ldots,m_{i-1},m_i-1,m_{i+1},\ldots\rangle , \nn
\eea
reads~:
\bea
&& \kappa_i |m\rangle=e^{i\pi m_i/k} |m\rangle,\nn\\
&&a_i^+ |m\rangle=e^{-i\pi(m_1+\ldots +m_{i-1})/k}
\sqrt{2\sin(\pi(m_i+1)/k)\sin(\pi/(2k))\over \sin^2(\pi/k)}
|m_i+1\rangle,  \label{action} \\
&&a_i^- |m\rangle=
e^{i\pi(m_1+\ldots +m_{i-1})/k}
\sqrt{2\sin(\pi m_i/k)\sin(\pi/(2k))\over \sin^2(\pi/k)}
 |m_i-1\rangle . \nn
\eea

It is also interesting to consider the decomposition of the
above Fock space representation, irreducible with respect to
$U_{q}[osp(1/2n)]$, according to the quantum subalgebra
$U_{q}[gl(n)]$. From~(\ref{eij}) we can in fact directly deduce
the representatives $\pi(e_{ij})$ of all Cartan-Weyl generators of
$U_{q}[gl(n)]$ in the Fock space~:
\bea
&&\pi(e_{ij})=-\cos(\pi/(2k))\kappa_ja^+_ja^-_i,\qquad{\rm for\ }i<j,
\\
&&\pi(e_{ij})=-\cos(\pi/(2k)) a_j^+ a_i^-
\kappa_i^{-1},\qquad{\rm for\ }i>j.
\eea
Then, the actual matrix elements follow from~(\ref{action}).
Note that the subspace spanned by vectors
$|m_1,\ldots,m_n\rangle$ with $m_1+\cdots+m_n=m$ ($m$ constant)
forms in fact a submodule for the $U_{q}[gl(n)]$ action.
Since the matrix elements~(\ref{action}) of $a_i^-$, resp.\
$a_i^+$, are nonzero for $m_i\ne 0$, resp.\ $m_i\ne k$, one
deduces that this $U_{q}[gl(n)]$ module is also irreducible. Thus, the
$k^n$-dimensional irreducible $U_{q}[osp(1/2n)]$ Fock
representation splits into $nk-n+1$ irreducible $U_{q}[gl(n)]$
representation $(m)$ characterized by a ``total number'' $m$ taking
values from $0$ up to $n(k-1)$. The dimension of this irreducible
representation $(m)$ is equal to the coefficient of $x^m$ in the
expansion of
\beq
(1+x+\cdots+x^{k-1})^n = \left({1-x^k\over 1-x}\right)^n,
\eeq
or, more explicitly,
\beq
\dim(m) = \sum_{j_0,j_1,\ldots,j_{k-1}} {n!\over j_0!j_1!\cdots
j_{k-1}!},
\eeq
where the $j_i$ assume all nonnegative integer values such that
$j_0+j_1 + \cdots + j_{k-1}=n$ and $j_1+2j_2 + \cdots +
(k-1)j_{k-1}=m$. Of course, the representations $(m)$ obtained
by means of deformed Bose operators are only a small part of the so-called
type~1 representations of $U_q[gl(n)]$ (see, e.g., Chapter~11 of
Ref.~\cite{Chari}, and references therein).

\section*{Acknowledgements}

We have pleasure in thanking Prof.\ R.\ Jagannathan and Dr.\
N.I.\ Stoilova for stimulating discussions.
One of us (TDP) is thankful to Prof.\ G.\ Vanden Berghe for the kind
hospitality at the Department of Applied Mathematics and Computer Science,
University of Ghent.

This research was supported by a grant awarded by the Belgian
National Fund for Scientific Research and by the contract
$\Phi$-416 of the Committee of Science of Bulgaria.


\begin{thebibliography}{99}
\bibitem{Palev78}
Ganchev A Ch and Palev T D 1978 {\it Preprint JINR} P2-11941;
1980 {\it J.\ Math.\ Phys.}
         {\bf 21} 797
\bibitem{Ohnuki}
Ohnuki Y and Kamefuchi S 1982 {\it Quantum Field Theory
and Parastatistics} (Univ.\ of Tokyo Press, Springer-Verlag,
Berlin)
\bibitem{Palev82}
 Palev T D 1982  {\it J.\ Math.\ Phys.} {\bf 23} 1100
\bibitem{Jag}
Jagannathan R and Vasudevan R 1984 {\it J.\ Math.\ Phys.} {\bf
25} 2294
\bibitem{Blank}
Blank J and  Havlicek M 1986  {\it J.\ Math.\ Phys.}
           {\bf 27} 2823
\bibitem{Biswas}
Biswas S N and Soni S K 1988 {\it J.\ Math.\ Phys.}
          {\bf 29} 16
\bibitem{FlatoF}
Flato M and Fronsdal C 1989 {\it Journ.\ of
                       Geometry and Physics} {\bf 6} 293
\bibitem{Okubo}
Okubo S 1994  {\it J.\ Math.\ Phys.} {\bf 35} 2785
\bibitem{Hughes}
Hughes J W B 1981 {\it J.\ Math.\ Phys.} {\bf 22} 245
\bibitem{Greenberg87}
Greenberg O W and Mohapatra R N 1987 {\it Phys.\
         Rev.\ Lett.} {\bf 59} 2507
\bibitem{Floreanini90}
Floreanini R and Vinet L 1990 {\it J.\ Phys.\ A~: Math.\
         Gen.} {\bf 23} L1019
\bibitem{Celeghini}
Celeghini E, Palev T D and Tarlini M 1990 {\it Preprint}
         YITP/K-865 Kyoto and
     1991  {\it Mod.\ Phys.\ Lett.\ B} {\bf 5} 187
\bibitem{Okada}
Odaka K, Kishi T and Kamefuchi S 1991
	 {\it J.\ Phys.\ A~: Math.\ Gen.} {\bf 24} L591
\bibitem{Beckers}
Beckers J and Debergh N 1991
          {\it J.\ Phys.\ A~: Math.\ Gen.} {\bf 247} L1277
\bibitem{Chaturvedi}
Chaturvedi S and Srinivasan V 1991 {\it Phys.\ Rev.\ A}
          {\bf 44} 8024
\bibitem{Palev93}
Palev T D 1993 {\it Lett.\ Math.\ Phys.} {\bf 28} 321
\bibitem{Krishna}
Krishna-Kumari M, Shanta P, Chaturvedi S and Srinivasan V
          1992 {\it Mod.\ Phys.\ Lett.\ A } {\bf 7 } 2593
\bibitem{Hadji}
Hadjiivanov L K 1993 {\it J.\ Math.\ Phys.} {\bf 34}
         5476
\bibitem{Bonatsos}
Bonatsos D and Daskaloyannis C 1993 {\it
          Phys.\ Lett.\ B} {\bf 307} 100 and the references therein
\bibitem{Flato93}
Flato M, Hadjiivanov L K and Todorov I T 1993
         {\it Found.\ Phys.} {\bf 23 } 571
\bibitem{Macfarlane93}
Macfarlane A J 1993  {\it Generalized Oscillator
          Systems and Their Parabosonic Interpretation 	}
 	  Preprint DAMPT 93-37
\bibitem{Palev93b}
Palev T D 1993 {\it J.\ Math.\ Phys.} {\bf 34} 4872
\bibitem{Palev93c}
Palev T D and Stoilova N I 1993
          {\it Lett.\ Math.\ Phys.} {\bf 28} 187
\bibitem{Quesne}
Quesne C 1994 {\it Phys.\ Lett.\ A} {\bf 193} 245
\bibitem{VdJJ}
Van der Jeugt J and Jagannathan R 1994 {\it Polynomial
deformations of $osp(1/2)$ and generalized parabosons} hep-th/9410145
\bibitem{Macfarlane94}
Macfarlane A J 1994 {\it  J.\ Math.\ Phys.} {\bf 35}
          1054
\bibitem{Cho}
Cho K H, Chaiho Rim, Soh D S and Park S U 1994
          {\it J.\ Phys.\ A~: Math.\ Gen.} {\bf 27} 2811
\bibitem{Chakrabarti}
Chakrabarti R and Jagannathan R 1994
         {\it J.\ Phys.\ A~: Math.\ Gen.} {\bf 27} L277
\bibitem{Palev94}
Palev T D 1994  {\it Lett.\ Math.\ Phys.} {\bf 31} 151
\bibitem{Green94}
Green H S 1994 {\it Austr.\ J.\ Phys.} {\bf 47} 109
\bibitem{Chaichian}
Chaichian M and Kulish P 1990 {\it Phys.\ Lett.\ B}
          {\bf 234} 72
\bibitem{Bracken}
Bracken A J, Gould M D and Zhang R B 1990, {Mod.\ Phys.\ Lett.\
A} {\bf 5} 331
\bibitem{Floreanini2}
Floreanini R, Spiridonov V P and Vinet L 1990
          {\it Phys.\ Lett.\ B} {\bf 242} 383
\bibitem{Floreanini3}
Floreanini R, Spiridonov V P and Vinet L 1991
          {\it Commun.\ Math.\ Phys.}{\bf 137} 149
\bibitem{Khoroshkin}
Khoroshkin S M  and Tolstoy  V N  1991 {\it Commun.\ Math.\ Phys.}
    {\bf 141} 599
\bibitem{dHoker}
d'Hoker E, Floreanini R and Vinet L 1991 {\it J.\ Math.\ Phys.}
{\bf 32} 1427
\bibitem{Green}
Green H S 1953 {\it Phys.\ Rev.} {\bf 90}  270
\bibitem{Omote}
Omote M, Ohnuki Y  and Kamefuchi S 1976 {\it Prog.\ Theor.\ Phys.}
 {\bf 56} 1948
\bibitem{Kac}
Kac V G 1978 {\it Lect.\ Notes in Math.} {\bf 626}  597
\bibitem{Greenberg}
Greenberg O W and Messiah A M 1965 {\it Phys.\ Rev.} {\bf 138B} 1155
\bibitem{Palev93A}
Palev T D 1993 {\it J.\ Phys.\ A~: Math.\ Gen.} {\bf 26} L1111
\bibitem{Palev91}
Palev T D and Tolstoy V N 1991 {\it Commun.\ Math.\ Phys.}
          {\bf 141} 549
\bibitem{Macfarlane}
Macfarlane A J, 1989 {\it J.\ Phys.\ A~: Math.\ Gen.} {\bf 22}  4581
\bibitem{Biedenharn}
Biedenharn L C, 1989 {\it J.\ Phys.\ A~: Math.\ Gen.} {\bf 22}  L873
\bibitem{Sun}
Sun C P and Fu H C, 1989 {\it J.\ Phys.\ A~: Math.\ Gen.} {\bf 22}  L983
\bibitem{Hayashi}
Hayashi T, 1990 {\it Commun.\ Math.\ Phys.} {\bf 127}  129
\bibitem{Pusz1}
Pusz W and Woronowicz S L 1989 {\it Rep.\ Math.\ Phys.} {\bf 27} 231
\bibitem{Pusz2}
Pusz W 1989 {\it Rep.\ Math.\ Phys.} {\bf 27} 349
\bibitem{Hadjiivanov}
Hadjiivanov L K, Paunov R R and Todorov I T 1992
          {\it J.\ Math.\ Phys.} {\bf 33} 1379
\bibitem{Jagannathan}
Jagannathan R, Sridhar R, Vasudevan R, Chaturvedi S,
          Krishnakumari M, Shanta P and Srinivasan V  1992
 {\it J.\ Phys.\ A~: Math.\ Gen.} {\bf 25} 6429
\bibitem{VdJ}
Van der Jeugt J 1993  {\it J.\ Phys.\ A~: Math.\ Gen.} {\bf 26} L405
\bibitem{Wess}
Wess J and Zumino B 1990 {\it Nucl.\ Phys.\ Proc.\ Suppl.\ B}
          {\bf 18} 302
\bibitem{Zumino}
Zumino B 1991 {Mod.\ Phys.\ Lett.\ A} {\bf 6} 1225
\bibitem{Furlan}
Furlan P, Hadjiivanov L K and Todorov I T 1992
          {\it J.\ Math.\ Phys.} {\bf 33} 4255
\bibitem{Flato}
Flato M, Hadjiivanov L K and Todorov I T 1993
          {\it Found.\ Phys.} {\bf 23} 571
\bibitem{Coon}
Coon D D, Yu S and Baker M 1972 {\it Phys.\ Rev.\ D}
          {\bf 15} 1429
\bibitem{Arik}
Arik M and Coon D D 1976 {\it J.\ Math.\ Phys.} {\bf 17} 524
\bibitem{Kuryshkin}
Kuryshkin V V 1980 {\it Ann.\ Fond.\ Louis de Broglie} {\bf 5}
          111
\bibitem{Jannussis}
Jannussis A, Brodimas G, Sourlas D and Zisis V 1981
          {\it Lett.\ Nuovo Cimento} {\bf 30} 123
\bibitem{Chari}
Chari V and Pressley A 1994 {\it A guide to Quantum Groups}
(Cambridge University Press, Cambridge)
\end{thebibliography}
\end{document}